\documentclass[11pt]{amsart}
\usepackage{amscd,amssymb}

\newcommand{\mC}{{\mathbb C}}

\newcommand{\mZ}{{\mathbb Z}}

\newcommand{\bO}{\Omega}
\newcommand{\bo}{\omega}

\newcommand{\mcD}{\mathcal D}

\newcommand{\mcH}{\mathcal H}

\newcommand{\mcL}{\mathcal L}

\newcommand{\mcN}{\mathcal N}
\newcommand{\mcO}{\mathcal O}
\newcommand{\mcP}{\mathcal P}

\newcommand{\mcR}{\mathcal R}
\newcommand{\mcS}{\mathcal S}
\newcommand{\mcT}{\mathcal T}

\newcommand{\mcZ}{\mathcal Z}

\newcommand{\fg}{\mathfrak g}

\newcommand{\fm}{\mathfrak m}

\newcommand{\ft}{\mathfrak t}

\newcommand{\ti}{\tilde}

\newcommand{\un}{\underline}
\newcommand{\ov}{\overline}

\title{On Shalika germs}

\author{David Kazhdan}
\address{Institute of Mathematics\\
The Hebrew University of Jerusalem\\
Givat-Ram, Jerusalem,  91904\\
Israel} \email{kazhdan@math.huji.ac.il}

\date{\today}
\begin{document}

\maketitle

\begin{abstract}
Let $G$ be a reductive group over a local field $F$ satisfying the assumptions of \cite{Deb1}, $G_{reg}\subset G$ the subset of regular elements. Let $T\subset G$ be a maximal torus. We write $T_{reg}=T\cap G_{reg}$. Let $dg ,dt$ be Haar measures on $G$ and $T$. They define an invariant measure $dg/dt$ on $G/T$. Let $\mcH$ be the space of complex valued locally constant functions on $G$ with compact support.
For any $f\in \mcH ,t\in T_{reg}$ we define $I_t(f)=\int _{G/T}f(\bar gt\bar g^{-1})dg/dt$.
Let $P$ be the set of conjugacy classes of unipotent
elements in $G$. For any $\bO \in P$ we fix an invariant measure
$\bo$ on $\bO$. As well known \cite {R} for any $f\in \mcH$ the integral
$$I_\bO (f)=\int _\bO f\bo$$ is absolutely convergent.
Shalika \cite{Sh} has shown that there exist functions $\ti j_\bO (t),\bO \in P$ on $T\cap G_{reg}$ such that
$$I_t(f)=\sum _{\bO \in P}\ti j_\bO (t)I_\bO(f) \eqno{(\star)}$$
for any $f\in \mcH ,t\in T$ {\it near} to $e$ where the notion of {\it near} depends on $f$.

 For any positive real number $r$ one defines an open $Ad$-invariant subset $G_r$ of $G$ and a subspace $\mcH _r$ as in \cite{Deb1}. In this paper I show that for any
 $f\in \mcH _r$ the equality $(\star)$ is true for all $t\in T_{reg}\cap G_r$.
 \end{abstract}

Let $F$ be a local non-archimedian field,  $\mcO$ be the ring of integers of $F, \fm \subset \mcO$ the maximal ideal, $k=\mcO / \fm$ the residue field, $\zeta$ a generator of $\fm, q=|k|$ and $val:F^\star \to \mZ$ the valuation such that $val(\zeta )=1$.

Let $G$ a  reductive  $F$-group, $G_{reg}\subset G$ the subset of regular elements. To simplify notations I write $G$ instead of $G(F)$ and fix a Haar
measure $dg$ on $G$.
For any $r\in \mZ _+$ we denote by $G_r$ the $Ad$-invariant
open subset as in \cite{Deb1}.

 Let $\fg$ be the Lie algebra of $G$ and $\fg _{nil}\subset \fg$ be the subset of elements $x$ of the form $x=n+z$ where $n$ is a topologically nilpotent element
of $[\fg ,\fg]$ and $z\in Z_\fg (\fm )$ where $Z_\fg$ is the center
of $\fg$.

 Let $\mcH =\mcS (G)$ be the space of locally constant $\mC$-valued compactly supported functions on $G$. Since we fixed a Haar measure on $G, \mcH$ has a natural algebra structure.
We denote by $\bar \mcH$ the space of locally constant $\mC$-valued compactly supported function on $\fg$.

I will freely use notations and results of \cite{Deb1}. In particular
we fix an map $\phi :\fg _n \to G$ as in \cite{Deb1}. (So
$\phi (x)=Id_n+x$ if $G=Gl_n$ and $\phi$ is the inverse of the quasi-logarithm map (see \cite[Appendix C]{BKV}) if $char (k)>>1$.)
and consider open $Ad$-invariant subsets $G_r\subset G, \fg _r\subset \fg$  as in \cite{Deb1}. Then $\phi (\fg _r)=G_r$ (see \cite[Lemma C.4]{BKV}). 

We fix a maximal torus $T$ in $G$ and  a Haar measure $dt$ on $T$. Let $\ft$ be the Lie algebra of $T, T_{reg}\subset T,\ft _{reg}\subset \ft$ the subsets of regular elements.

 For any regular $t\in T_{reg} ,  f\in \mcH$ we define
$$I_t(f)=\int _{G/T} f(g\phi (t)g^{-1})dg/dt.$$
Analogously for any regular $\bar t\in \ft_{reg} ,  \bar f\in \bar \mcH$ we define
$$I_{\bar t}(\bar f)=\int _{G/T} \bar f(g\bar tg^{-1})dg/dt.$$

It is clear that for any  $f\in \mcH ,\bar t\in \ft _{reg}\cap \fg _{nil}$ we have
$$I_t(f)=I_{\bar t}(\phi ^\star (f)),t=\phi (\bar t).$$

We define $\mcH _r\subset \mcH$ as the subspace of $f\in \mcH$
such that $ \pi (f)=0$
for all irreducible representations $\pi$ of $G$ of {\it depth}
$> r$. As follows from
\cite[Teorem 5.5.4]{Deb2} this subspace  coincides with the subspace $\mcH _r$ of
\cite{Deb1}.

Let $P$ be the set of conjugacy classes of nilpotent
elements in $G$. For any $\bO \in P$ we fix an invariant measure
$\bo$ on $\bO$. As well known \cite {R} for any $\bar f\in \bar \mcH$ the integral
$$I_\bO (\bar f)=\int _\bO \bar f\bo$$
is absolutely convergent and
$I_\bO (\bar f_\zeta  )=q^dI_\bO (\bar f)$ where $\bar f_\zeta (x)=f(\zeta ^2 x)$.

The following result is proven in \cite[Conjecture 2]{Deb1}.

{\bf Claim}. For any $f\in \mcH _r$ such that $I_\bO (\phi ^\star (f))=0$ for all $\bO \in P$ we have
$I_t(f)=0$ for all $t\in T_{reg}\cap G_r$.

 Shalika \cite{Sh} has shown that there exist functions $\ti j_\bO (t),\bO \in P$ on $T\cap G_{reg}$ such that
$$I_t(f)=\sum _{\bO \in P}\ti j_\bO (t)I_\bO (f)$$
for any $t\in T$ near to $e$. Let $j_\bO$ be the function
on a neighboorhood of $0$ in $\mcT$ given by $j_\bO (x)=\ti j_\bO (\phi (x))$.

Harish-Chandra \cite{HC} has shown that

 $$j_\bO (c^2x)=j_\bO (x)|c|^{dim (\mcN )-dim (\bO)}$$
for all $c\in F^\star$ if $x\in \mcT$ is near to $0$.

Using this homogenuity we can define functions $j_\bO (t)$ for all $t\in \mcT$.

{\bf Theorem}. For any $f\in \mcH _r$ the Shalika germ expansion
is true for all $g\in \mcT \cap \fg _r$.

{\bf Proof}. For any $\bO \in P$ we define $\mcH _r (\bO )$ to be the set of all $f\in \mcH _r$ such that  $I_{\bO '} (f)=0$ for all $\bO '\in P$ different from $\bO$. 

As follows from \cite[Theorem 4.1.4]{Deb1} the space
$\mcH _r$ is spanned by subspaces $\mcH _r(\bO ), \bO \in P$.
So it is sufficient for a  proof of  Theorem
to show that
$$q^dI_{\bar t}(\bar f)=I_{\zeta ^2\bar t}(\bar f),\bar f=\phi ^\star (f)$$
for all $f\in \mcH _r(\bO )$ all regular  $\bar t\in \ft \cap \fg _r$.

Let
$$\bar f=\phi ^\star (f), h=q^d \bar f-\bar f_\zeta $$
where  $d=dim (\bO)$.
Then $I_{\bO '}(h)=0$ for all $\bO ' \in P$ and it follows from the  Claim that $I_t(h)=0$
for all $t\in  \mcT _r$. So
$$q^dI_{\bar t}(\bar f)=I_{\zeta ^2\bar t}(\bar f),\bar f=\phi ^\star (f).$$

\end{document}